\documentclass[preprint,12pt]{elsarticle}

\usepackage{afterpage}
\usepackage{amssymb}
\usepackage{amsmath}
\usepackage{amsthm}
\usepackage{makecell}

\usepackage{subfigure}
\usepackage{algorithm}
\usepackage{algpseudocode}
\usepackage{listings}
\usepackage{xcolor}

\usepackage{bm}
\usepackage{geometry}
\usepackage{lineno}

\usepackage{booktabs} 

\newcommand{\norm}[1]{\left\lVert #1 \right\rVert}

\newtheorem{Lemma}{Lemma}
\newtheorem{Conjecture}{Conjecture}
\newtheorem{Theorem}{Theorem}

\newtheorem{example}{Example}

\journal{Applied Mathematics and Computation}

\begin{document}

\begin{frontmatter}

\title{Counterexamples to a Conjecture on First Derivative Bounds of Rational Bézier Curves}

\author{Mao Shi}
\affiliation{organization={School of Mathematics and Statistics, Shaanxi Normal University},
            city={Xi'an},
            country={China}}

\begin{abstract}
In this paper we present an explicit counterexample of degree $n=7$, which shows that the conjecture proposed by Li et al. \cite{Li2013} regarding the first derivative bounds for rational Bézier curves is generally false.
We further derive an explicit rational B\'ezier representation of the first derivative and propose a degree-elevation based computable upper bound for $\sup_{t\in[0,1]}\|\mathbf r'(t)\|$. 
The bound is valid for any finite elevation order and converges to the true supremum as the elevation degree tends to infinity.
An \emph{a priori} tolerance-driven rule is provided to determine a sufficient elevation degree, and the computational complexity of the proposed procedure is analyzed.
Numerical experiments validate the counterexample and demonstrate the accuracy and efficiency of the new upper bound across a range of degrees and weight patterns.

\end{abstract}

\begin{highlights}
\item A non-degenerate counterexample of degree n=7 disproves a conjectured first-derivative bound for rational Bézier curves.
\item An explicit rational Bézier representation of $r'(t)$ is derived via Bernstein-product and degree-elevation formulas.
\item A computable degree-elevation upper bound for $\sup_{t\in[0,1]}\|\mathbf r'(t)\|$ is established for any finite elevation order.
\item The proposed bound converges to the true supremum as the elevation degree increases, with a quantified convergence rate.
\item A tolerance-driven criterion is provided to select a sufficient elevation order, supported by numerical experiments.
\end{highlights}

\begin{keyword}
Rational Bézier curve \sep Derivative bounds \sep Counterexamples \sep Supremum analysis 
\end{keyword}

\end{frontmatter}

\section{Introduction}
\label{sec1}

Rational Bézier curves serve as essential tools in geometric modeling applications, where accurate derivative estimation plays a vital role in computational procedures. Maintaining control over derivative magnitudes ensures numerical stability during curve evaluation, facilitates adaptive subdivision approaches, and improves collision detection mechanisms \cite{Farin2002,Sederberg1987,Floater1992,Wang1995}.

The standard representation of an $n$th-degree rational Bézier curve takes the form \cite{Farin2002}
\begin{equation}\label{eq:rational_bezier}
  \mathbf{r}\left( t \right) = \frac{\mathbf{p}(t)}{\omega(t)} = \frac{\sum\limits_{i=0}^n{\omega _i\mathbf{r}_iB_{i}^{n}\left( t \right)}}{\sum\limits_{i=0}^n{\omega _iB_{i}^{n}\left( t \right)}},\ t\in \left[0,1\right],
\end{equation}
where $B_i^n(t)$ represent the Bernstein basis polynomials, $\mathbf{r}_i$ denote control points in $\mathbb{R}^d$, and $0 <\omega_i < +\infty$ indicate the corresponding weights.

Let $M$ denote the weight variation measure
\begin{equation}\label{eq:weight_parameter}
  M = \max \left\{ \max_{0 \le i \le n-1} \frac{\omega_i}{\omega_{i+1}}, \max_{0 \le i \le n-1} \frac{\omega_{i+1}}{\omega_i} \right\},
\end{equation}
Zhang and Ma demonstrated \cite{Zhang2006} the bounding relationship
\begin{equation} \label{eq:zhang_ma_bound}
\norm{\mathbf{r}'(t)} \leqslant nM \max_{0 \le i \le n-1} \{\norm{\mathbf{r}_{i+1} - \mathbf{r}_i}\},
\end{equation}
valid for degrees $n=2,3$, where $\|\cdot\|$ denotes the Euclidean norm.

Li et al. \cite{Li2013} extended this observation by proposing the following conjecture for general degrees:
\begin{Conjecture}\label{conj:Li}
For rational Bézier curves of arbitrary degree $n \in \mathbb{N}$, the inequality \eqref{eq:zhang_ma_bound} holds universally.
\end{Conjecture}
This conjecture was further supported by Zhang and Jiang \cite{Zhang2016}. 

In this paper, we show that Conjecture~\ref{conj:Li} fails at degree $n=7$  in Section \ref{sec:counterexample}, hence it is not valid in general.  We also investigate a convergent degree-elevation upper bound for
\(\sup_{t\in[0,1]}\norm{\mathbf{r}'(t)}\) in Section~\ref{sec:supremum}. Finally, we give more examples in Section \ref{Numerical Validation}.

\section{Counterexample}
\label{sec:counterexample}

\begin{figure}[htbp]
    \centering
    \begin{minipage}[b]{0.5\textwidth}
        \centering
        \includegraphics[width=\textwidth]{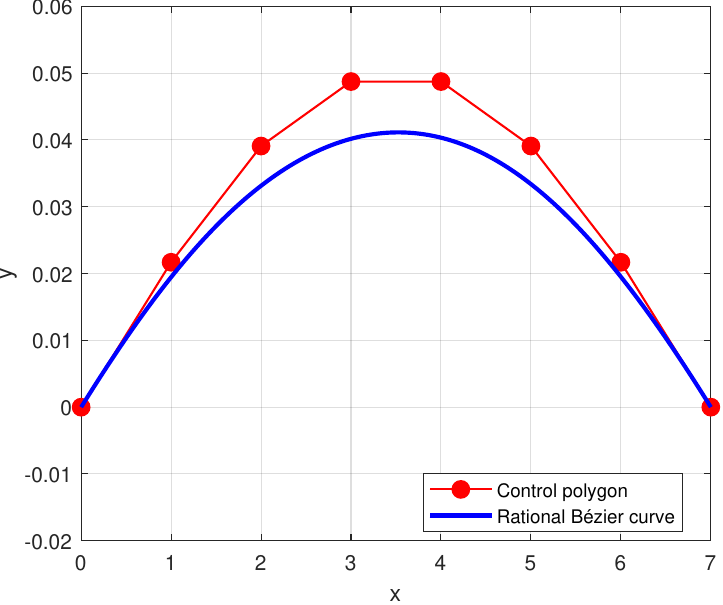}
        \vspace{2mm}
        \footnotesize (a) 
    \end{minipage}
    \hfill
    \begin{minipage}[b]{0.5\textwidth}
        \centering
        \includegraphics[width=\textwidth]{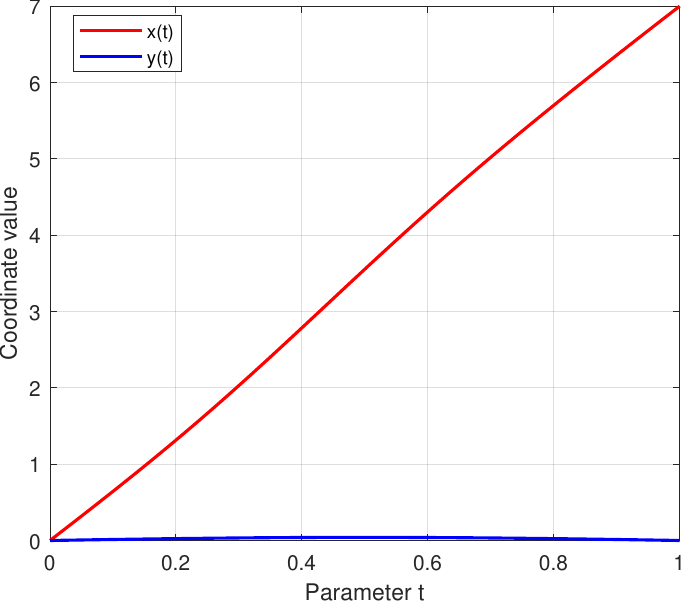}
        \vspace{2mm}
        \footnotesize (b) 
    \end{minipage}
    \hfill
    \begin{minipage}[b]{0.5\textwidth}
        \centering
        \includegraphics[width=\textwidth]{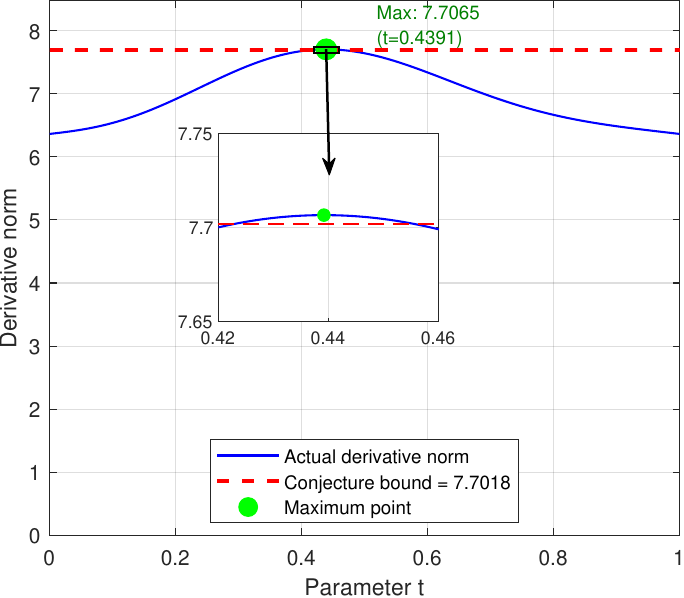}
        \vspace{2mm}
        \footnotesize (c) 
    \end{minipage}
    
    \vspace{3mm}
    \caption{Schematic of Example 1. (a) 7th-degree rational Bézier curve. (b) X-Y components of the rational Bézier curve. (c) Comparison of the first-derivative norm curve with the conjectured upper bound.}
    \label{Fig:CounterexampleB}
\end{figure}

\begin{example}
   Consider a rational Bézier curve of degree $ n = 7$ with weights

\[
\boldsymbol{\omega} = (1.3310, 1.2100, 1.1000, 1.0000, 1.1000, 1.2100, 1.3310, 1.4641)
\]
and the control points
\[
\begin{split}
\mathbf{r} &= ((0.00,0.0000), (1.00,0.0217), (2.00,0.0391), (3.00,0.0487), \\
           &\quad \ \ \  (4.00,0.0487), (5.00,0.0391), (6.00,0.0217), (7.00,0.0000)).
\end{split}
\]

The norm of the corresponding first-order derivative can be expressed as

\[
\norm{\mathbf{r}'(t)} = \frac{\norm{\mathbf{p}'(t)\omega(t)-\mathbf{p}(t)\omega'(t)}}{\omega(t)^2}=\frac{\sqrt{Q(t)}}{P(t)^4},
\]
where
\[
P(t) = - 13.31+ 8.47 t- 2.31 t^{2}+ 0.35 t^{3}- 66.85 t^{4}+ 119.91 t^{5}- 79.94 t^{6}+ 19.039 t^{7},
\]
and
\begin{tiny}
\[
\begin{split}
Q(t) &=- 2.775543946\times 10^{6} t+ 2.644930730\times 10^{9} t^{15}- 2.085446174\times 10^{9} t^{16}+ 9.218380764\times 10^{8} t^{17}- 1.436345689\times 10^{8} t^{18} \\
     &- 1.146974444\times 10^{8} t^{19}+ 1.028576611\times 10^{8} t^{20}- 4.250234212\times 10^{7} t^{21}+ 1.039418663\times 10^{7} t^{22}- 1.453809244\times 10^{6} t^{23} \\
     &+ 91515.669 t^{24}+ 2.901398784\times 10^{6} t^{2} + 3.016301925\times 10^{7} t^{3}- 1.476249284\times 10^{8} t^{4}+ 3.231094005\times 10^{8} t^{5}- 2.192757049\times 10^{8} t^{6} \\
     &- 1.070822330\times 10^{9} t^{7}+ 4.627429383\times 10^{9} t^{8}- 1.003934933\times 10^{10} t^{9}+ 1.441684279\times 10^{10} t^{10}- 1.451404144\times 10^{10} t^{11} \\
     &+ 9.985885877\times 10^{9} t^{12}- 3.690810465\times 10^{9} t^{13}- 9.934198245\times 10^{8} t^{14}+ 1.271531954\times 10^{6}
\end{split}
\]
\end{tiny}
Through numerical optimization, such as the function 'fminbnd' in Matlab(R), we find
\[
\max_{t\in[0,1]}\norm{\mathbf{r}'(t)} = 7.7065,
\]
which occurs at $t = 0.4391$. 

By substituting the relevant weights $\boldsymbol{\omega}$ into Equation \eqref{eq:weight_parameter}, we obtain that $M = 1.100$, which implies that the right-hand side of the conjecture is
 $$nM \max_{0 \le i \le n-1} \{\norm{\mathbf{r}_{i+1} - \mathbf{r}_i}\}=7.7018.$$ 
 
 Since $7.7065 > 7.7018$, we have constructed a valid counterexample to refute the conjecture (See Fig.~\ref{Fig:CounterexampleB}).
\end{example}

\section{Supremum Analysis of $\norm{\mathbf{r}'(t)}$}
\label{sec:supremum}

Building on established methodologies for derivative representations of rational B\'ezier curves \cite{Shi2025}, we derive explicit expressions for the first derivative and investigate the quantity
\(\sup_{t\in[0,1]}\norm{\mathbf r'(t)}\).
The same framework extends naturally to derivatives of higher order.

\subsection{First Derivative Representation}

We begin by deriving an explicit representation for the first derivative of a rational Bézier curve.
\begin{Theorem}
\label{thm:derivative_explicit}
The first derivative of a rational Bézier curve permits the representation:
\begin{equation}\label{eq:derivative_rep}
\mathbf{r}'(t) = \frac{n\sum_{i=0}^{2n} \widehat{\mathbf{R}}_{i,2n}^{[1]} B_i^{2n}(t)}{\sum_{i=0}^{2n} \omega_{i,2n}^{[1]} B_i^{2n}(t)} = \frac{\sum_{i=0}^{2n} \mathbf{Q}_{i,2n}^{[1]} \omega_{i,2n}^{[1]} B_i^{2n}(t)}{\sum_{i=0}^{2n} \omega_{i,2n}^{[1]} B_i^{2n}(t)},
\end{equation}
with constituent parameters determined through the following computational procedures:

\begin{enumerate}
\item \textbf{Weight coefficients} $\omega^{[1]}_{i,2n}$ for $i = 0, \ldots, 2n$:
\begin{equation}\label{eq:weight_coeff}
\omega^{[1]}_{i,2n} = \sum_{j=\max(0,i-n)}^{\min(i,n)} \frac{\binom{n}{j}\binom{n}{i-j}}{\binom{2n}{i}} \omega_j \omega_{i-j}.
\end{equation}

\item \textbf{Intermediate control points} $\mathbf{R}^{[1]}_{j,2n}$ for $j = 0, \ldots, 2n-1$:
\begin{equation}\label{eq:intermediate_R}
\mathbf{R}^{[1]}_{j,2n} = \sum_{h=\max(0,j-n)}^{\min(n-1,j)} \frac{\binom{n-1}{h}\binom{n}{j-h}}{\binom{2n-1}{j}} \left[\omega_{h+1}\omega_{j-h}(\mathbf{r}_{h+1} - \mathbf{r}_{j-h}) + \omega_h \omega_{j-h}(\mathbf{r}_{j-h} - \mathbf{r}_h)\right].
\end{equation}

\item \textbf{Elevated control points} $\widehat{\mathbf{R}}^{[1]}_{i,2n}$ for $i = 0, \ldots, 2n$:
\begin{equation}\label{eq:Rhat}
\widehat{\mathbf{R}}^{[1]}_{i,2n} = \sum_{j=\max(0,i-1)}^{\min(2n-1,i)} \frac{\binom{2n-1}{j}}{\binom{2n}{i}} \mathbf{R}^{[1]}_{j,2n}.
\end{equation}

\item \textbf{Final control points} $\mathbf{Q}^{[1]}_{i,2n} $ for $i = 0, \ldots, 2n$:
\begin{equation}\label{eq:final_Q}
\mathbf{Q}^{[1]}_{i,2n} = \frac{n \widehat{\mathbf{R}}^{[1]}_{i,2n}}{\omega^{[1]}_{i,2n}}.
\end{equation}
\end{enumerate}
\end{Theorem}
\begin{proof}
From \eqref{eq:rational_bezier}, it yields
\begin{equation}\label{eq:quotient}
\mathbf{r}'(t) = \frac{\mathbf{p}'(t)\omega(t) - \mathbf{p}(t)\omega'(t)}{\omega(t)^2}.
\end{equation}
 Using the product formula for Bernstein polynomials \cite{Farouki1988}, we express $\omega(t)^2$ as
\[
\omega(t)^2 = \left(\sum_{i=0}^n \omega_i B^n_i(t)\right)^2 = \sum_{i=0}^{2n} \omega^{[1]}_{i,2n} B^{2n}_i(t),
\]
where $\omega^{[1]}_{i,2n}$ is given by \eqref{eq:weight_coeff}. 

Similarly, the numerator $\mathbf{p}'(t)\omega(t) - \mathbf{p}(t)\omega'(t)$ can be expressed as a polynomial of degree $2n-1$ in Bernstein form with coefficients $\mathbf{R}^{[1]}_{j,2n}$ given by \eqref{eq:intermediate_R}. After degree elevation from degree $2n-1$ to degree $2n$ using the standard elevation formula, we obtain $\widehat{\mathbf{R}}^{[1]}_{i,2n}$ in \eqref{eq:Rhat}.

Combining these results yields the representation \eqref{eq:derivative_rep}, completing the proof.
\end{proof}
This computational framework, consistent with established approaches \cite{Piegl1997,piegl1997a}, guarantees numerical reliability. 
Beyond its computational validity, the explicit formulation also admits a clear geometric interpretation. 

 The new weights $\omega_{i,2n}^{[1]}$ are derived from the product of two copies of the original weight function, reflecting the quadratic dependence on the original weights induced by differentiation of a rational function. They govern the rational parameterization of the derivative representation. The intermediate control vectors $\mathbf{R}_{j,2n}^{[1]}$ are linear combinations of scaled difference vectors $(\mathbf{r}_{k+1} - \mathbf{r}_k)$ from the original control polygon. Each vector aggregates the contributions of various polygon edges to the tangent direction at a specific parameter interval. After degree elevation to match the denominator, the final points $\mathbf{Q}_{i,2n}^{[1]}$ serve as the control points for this $2n$-degree rational representation. Geometrically, the rational function defined by $\left( \mathbf{Q}_{i,2n}^{[1]}, \omega_{i,2n}^{[1]} \right)$ does not trace a curve in the standard point space. Instead, its evaluated value at parameter $t$ is precisely the \emph{tangent vector} $\mathbf{r}'(t)$ of the original curve. Thus, these points provide a \emph{control structure} for the tangent vector field along the curve.

\subsection{Degree Elevation and Upper Bounds}

To establish rigorous upper bounds, we need the following two lemmas.

\begin{Lemma}\label{lem:elevation}
For any non-negative integer $e$, the first derivative admits the elevated representation:
\begin{equation}\label{eq:elevated}
\mathbf{r}'(t) = \frac{n \sum_{i=0}^{2n+e} \widehat{\mathbf{R}}^{[1]}_{i,2n+e} B^{2n+e}_i(t)}{\sum_{i=0}^{2n+e} \omega^{[1]}_{i,2n+e} B^{2n+e}_i(t)},
\end{equation}
where
\begin{equation}\label{eq:elevated_R}
\widehat{\mathbf{R}}^{[1]}_{i,2n+e} = \frac{1}{\binom{2n+e}{i}} \sum_{j=\max(0,i-e)}^{\min(2n,i)} \binom{2n}{j}\binom{e}{i-j} \widehat{\mathbf{R}}^{[1]}_{j,2n},
\end{equation}
and
\begin{equation}\label{eq:elevated_omega}
\omega^{[1]}_{i,2n+e} = \frac{1}{\binom{2n+e}{i}} \sum_{j=\max(0,i-e)}^{\min(2n,i)} \binom{2n}{j}\binom{e}{i-j} \omega^{[1]}_{j,2n}.
\end{equation}
\end{Lemma}

\begin{proof}
The degree elevation formula for Bernstein polynomials states that \cite{Farouki1988}:
\[
B^m_j(t) = \sum_{i=j}^{j+e} \frac{\binom{m}{j}\binom{e}{i-j}}{\binom{m+e}{i}} B^{m+e}_i(t).
\]
Applying this to both the numerator and denominator of \eqref{eq:derivative_rep} and collecting terms yields \eqref{eq:elevated_R} and \eqref{eq:elevated_omega}.
\end{proof}

\begin{Lemma}
For $\textbf{a}_k \in \mathbb{R}^d$, $b_k > 0$, and $c_k > 0$ with $1 \leq k \leq N$,
\begin{equation}\label{eq:fundamental}
\norm{ \frac{\sum_{k=1}^N \textbf{{a}}_k c_k}{\sum_{k=1}^N b_k c_k} } \leq \max_{k=1,\ldots,N} \frac{\norm{\textbf{a}_k}}{b_k}.
\end{equation}
\end{Lemma}
\begin{proof}
By virtue of the triangle inequality property of the norm and the following inequality \cite{Mitrinovic1970,Kuang1989}
\[
\frac{\sum a_k c_k}{\sum b_k c_k} \leqslant \max_k \frac{a_k}{b_k},
\]
where $a_k \in \mathbb{R}$. We can establish  \eqref{eq:fundamental}.
\end{proof}
We now establish the main upper bound result.

\begin{Theorem}\label{thm:bound}
The first derivative magnitude satisfies:
\begin{equation}\label{eq:bound}
\norm{\mathbf{r}'(t)} \leq \max_{i=0,\ldots,2n+e} \frac{n \norm{ \sum_{j=\max(0,i-e)}^{\min(2n,i)} \binom{2n}{j}\binom{e}{i-j} \widehat{\mathbf{R}}^{[1]}_{j,2n} }}{\sum_{j=\max(0,i-e)}^{\min(2n,i)} \binom{2n}{j}\binom{e}{i-j} \omega^{[1]}_{j,2n}}.
\end{equation}
Moreover, as $e \to +\infty$, this bound converges to $\sup\limits_{t \in [0,1]} \norm{\mathbf{r}'(t)}$.
\end{Theorem}

\begin{proof}
From \eqref{eq:elevated}, we have:
\[
\mathbf{r}'(t) = \frac{\sum_{i=0}^{2n+e} \left(n\widehat{\mathbf{R}}^{[1]}_{i,2n+e}\right) B^{2n+e}_i(t)}{\sum_{i=0}^{2n+e} \omega^{[1]}_{i,2n+e} B^{2n+e}_i(t)}.
\]
Since $B^{2n+e}_i(t) \geq 0$ for $t \in [0,1]$ and $\omega^{[1]}_{i,2n+e} > 0$, applying \eqref{eq:fundamental} with $a_i = n\widehat{\mathbf{R}}^{[1]}_{i,2n+e}$, $b_i = \omega^{[1]}_{i,2n+e}$, and $c_i = B^{2n+e}_i(t)$ yields:
\[
\norm{\mathbf{r}'(t)} \leq \max_{i=0,\ldots,2n+e} \frac{n\norm{\widehat{\mathbf{R}}^{[1]}_{i,2n+e}}}{\omega^{[1]}_{i,2n+e}}.
\]
Substituting \eqref{eq:elevated_R} and \eqref{eq:elevated_omega} gives \eqref{eq:bound}.

Furthermore, the control polygon defined by points $\mathbf{Q}^{[1]}_{i,2n+e} = n\widehat{\mathbf{R}}^{[1]}_{i,2n+e}/\omega^{[1]}_{i,2n+e}$ converges uniformly to the curve $r'(t)$ as $e \to +\infty$ \cite{Farin2002}. Consequently,
\[
\lim_{e \to +\infty} \max_{i=0,\ldots,2n+e} \norm{\mathbf{Q}^{[1]}_{i,2n+e}} = \sup_{t \in [0,1]} \norm{\mathbf{r}'(t)}.
\]
This completes the proof.
\end{proof}

Starting from the original rational Bézier data, the proposed method consists of a preprocessing stage
and a convolution-based degree-elevation stage.
The preprocessing, which constructs the intermediate weights and control points
$\omega^{[1]}_{j,2n}$ and $\widehat{\mathbf R}^{[1]}_{j,2n}$, requires $O(n^2 d)$ time.
The degree-elevation stage can be implemented as a direct convolution, involving exactly $(2n+1)(e+1)$
accumulation steps, each with $O(d)$ vector operations, yielding $O((2n+1)(e+1)d)$ time.
The additional costs for normalization and norm evaluation are linear in $(2n+e)d$ and do not affect
the leading order.
Overall, the time complexity is $O(n^2 d + (2n+1)(e+1)d)$, while the space complexity is
$O((n+e)d + n+e)$.

\begin{Theorem}
The degree-elevation upper bound for the first derivative converges to \(\sup_{t\in[0,1]}\norm{\mathbf r'(t)}\) at rate \(O\!\left(1/(2n+e)\right)\). 
\end{Theorem}
\begin{proof}
For a rational Bézier curve, its first derivative $\mathbf{r}'(t)$ can be represented in the canonical form of rational Bézier curves \eqref{eq:derivative_rep}; consequently, the derivative curve satisfies the same convergence rate for degree elevation as the original curve \cite{prautzsch1994}, which is given by $O(1/(2n+e))$.
\end{proof}

\subsection{Determining the Minimal Elevation Degree for a Given Tolerance}

To apply the degree elevation method in practice, it is crucial to determine the minimal elevation degree $e$ such that the computed upper bound lies within a user-specified tolerance $\varepsilon$ of the true supremum $\sup_{t\in[0,1]}\norm{\mathbf{r}'(t)}$. In this section, we extend the sharp quantitative bounds of Nairn, Peters and Lutterkort \cite{Nairn1999} from polynomial to rational Bézier curves, and use the resulting \emph{a priori} error bound to solve this tolerance satisfaction problem.

\begin{Theorem} \label{Vec3.1}
  Let $\mathbf{p}(t)$ be an $n$th-degree B\'ezier curve in $\mathbb{R}^d$ with control points \(\mathbf b_0,\dots,\mathbf b_n\) and $\bm{\ell}(t)$ its control polygon. Then
  \[
  \|\mathbf{p}-\bm{\ell}\| \le \mathbb{N}_{\infty}(n) \cdot \max_{0<i<n} \|\Delta_2\mathbf{b}_i\|,
  \]
  where $\Delta_2\mathbf{b}_i = \mathbf{b}_{i-1}-2\mathbf{b}_i+\mathbf{b}_{i+1}$, and the constant is
  \[
  \mathbb{N}_{\infty}(n)=\frac{\lfloor n/2\rfloor\lceil n/2\rceil}{2n}.
  \]
  \end{Theorem}
\begin{proof}
 The proof follows the vector-valued extension of Theorem 3.1 in \cite{Nairn1999}.
\end{proof}
\begin{Theorem} \label{Vec7.1}
  Let $\mathbf{p}(t)$ be a B\'ezier curve of degree $n$ in $\mathbb{R}^d$ with control point sequence $\mathbf{b}^{(n)}_0,\dots,\mathbf{b}^{(n)}_n$, and let
  \[
  M_n:=\max_{0<i<n}\|\Delta_2\mathbf{b}^{(n)}_i\|.
  \]
  Raise the degree of $\mathbf{p}$ to $m>n$, obtaining new control points $\mathbf{b}^{(m)}_0,\dots,\mathbf{b}^{(m)}_m$ and the corresponding control polygon $\bm{\ell}^{(m)}$. Then
  \[
  \|\mathbf{p}-\bm{\ell}^{(m)}\|\le \frac{n(n-1)}{m(m-1)}\,\mathbb{N}_{\infty}(m)\,M_n.
  \]
\end{Theorem}

\begin{proof}
  Consider successive single-step degree elevations. When raising from degree $k$ to $k+1$, the maximum second difference satisfies\cite{Nairn1999}
  \[
  \|\Delta_2\mathbf{b}^{(k+1)}\|\le\frac{k-1}{k+1}M_k.
  \]
    Applying the vector-valued version of Theorem \ref{Vec3.1}  to the curve of degree $k+1$ gives
  \[
  \|\mathbf{p}-\bm{\ell}^{(k+1)}\|\le\mathbb{N}_{\infty}(k+1)M_{k+1}.
  \]

  Combining these two estimates, we obtain the recurrence
  \[
  \|\mathbf{p}-\bm{\ell}^{(k+1)}\|\le\frac{k-1}{k+1}\cdot\frac{\mathbb{N}_{\infty}(k+1)}{\mathbb{N}_{\infty}(k)}\cdot\mathbb{N}_{\infty}(k)M_k.
  \]

  Iterating from $k=n$ to $k=m-1$ yields
  \[
  \|\mathbf{p}-\bm{\ell}^{(m)}\|\le
  \Bigl(\prod_{k=n}^{m-1}\frac{k-1}{k+1}\Bigr)
  \Bigl(\prod_{k=n}^{m-1}\frac{\mathbb{N}_{\infty}(k+1)}{\mathbb{N}_{\infty}(k)}\Bigr)
  \mathbb{N}_{\infty}(n)M_n.
  \]

  The first product telescopes:
  \[
  \prod_{k=n}^{m-1}\frac{k-1}{k+1}=\frac{(n-1)n}{(m-1)m}.
  \]
  The second product also telescopes:
  \[
  \prod_{k=n}^{m-1}\frac{\mathbb{N}_{\infty}(k+1)}{\mathbb{N}_{\infty}(k)}=\frac{\mathbb{N}_{\infty}(m)}{\mathbb{N}_{\infty}(n)}.
  \]
  Substituting and canceling $\mathbb{N}_{\infty}(n)$ gives
  \[
  \|\mathbf{p}-\bm{\ell}^{(m)}\|\le\frac{n(n-1)}{m(m-1)}\,\mathbb{N}_{\infty}(m)\,M_n.
  \]
  Sharpness follows from the case where all second-difference vectors are equal.
\end{proof}

\begin{Theorem}
  Let $\mathbf{r}(t)$ be a rational B\'ezier curve of degree $n$ with positive weights $w_i$ and control points $\mathbf{r}_i\in\mathbb{R}^d$. Raise its degree to $m>n$, yielding a new control polygon $\tilde{\bm{\ell}}^{(m)}$. Then
  \[
  \|\mathbf{r}-\tilde{\bm{\ell}}^{(m)}\|\le\frac{n(n-1)}{m(m-1)}\mathbb{N}_{\infty}(m)\,
  \frac{\displaystyle\max_{0<i<n}\Bigl(\|\Delta_2(\omega\mathbf{r})_i\|+\|\Delta_2 \omega\|_\infty\max_{0 \le j \le n}\|\mathbf{r}_j\|\Bigr)}{\displaystyle\min_{0\le i\le n}\omega_i},
  \]
  where  $\|\Delta_2 \omega\|_\infty=\max\limits_{1\le i\le n-1}|\Delta_2 \omega_i|$.
\end{Theorem}

\begin{proof}[proof]
  In homogeneous space $\mathbb{R}^{d+1}$, $\mathbf{r}^h(t) = \left(\mathbf{p}(t), \omega (t)\right)$ is a polynomial curve. Apply Lemma \ref{Vec7.1} to $\mathbf{r}^h$:
  \[
  \|\mathbf{r}^h-\tilde{\bm{\ell}}^{h,(m)}\|_\infty\le\frac{n(n-1)}{m(m-1)}\mathbb{N}_{\infty}(m)\max_i\|\Delta_2\mathbf{r}_i^h\|.
  \]
  
  Using the Euclidean norm in $\mathbb{R}^{n+1}$ and the triangle inequality, we get
\[
\|\Delta_2\mathbf{r}_i^h\|
\le \max_{0<i<n}\|\Delta_2(\omega\mathbf{r})_i\| + \max_{0<i<n}|\Delta_2 \omega_i|
= \max_{0<i<n}\|\Delta_2(\omega\mathbf{r})_i\| + \|\Delta_2 \omega\|_{\infty}.
\]
  
  Now we relate the distance in the original space $\mathbb{R}^d$ to that in homogeneous space.
  Set
\[
\mathbf{p}(t)=\mathbf{r}^h(t)_{\text{first }d}, \quad  \omega(t)=\mathbf{r}^h(t)_{\text{last}},
\qquad
\tilde{\mathbf{p}}(t)=\tilde{\bm{\ell}}^{h,m}(t)_{\text{first }d},\quad \tilde{\omega}(t)= \tilde{\bm{\ell}}^{h,m}(t)_{\text{last}}.
\]
  We have
  \[
  \Bigl\|\frac{\mathbf{p}(t)}{\omega(t)}-\frac{\tilde{\mathbf{p}}(t)}{\tilde{\omega}(t)}\Bigr\|
  \le\frac{\|\mathbf{p}(t)-{\tilde{\mathbf{p}}(t)}\|}{\omega(t)}+\frac{\|{\tilde{\mathbf{p}}(t)}\|\,|\omega(t)-\tilde{\omega}(t)|}{\omega(t)\tilde{\omega}(t)}.
  \]
  Using the homogeneous space bounds:
  \[
  \|\mathbf{p}(t)-\tilde{\mathbf{p}}(t)\|\le\frac{n(n-1)}{m(m-1)}\mathbb{N}_{\infty}(m)\max_i\|\Delta_2(\omega\mathbf{r})_i\|_2,
  \]
  \[
  |\omega(t)-\tilde{\omega}(t)|\le\frac{n(n-1)}{m(m-1)}\mathbb{N}_{\infty}(m)\max_i|\Delta_2 \omega_i|.
  \]
  The denominators $x_{d+1},y_{d+1}$ lie between $\min w_i$ and $\max w_i$; a conservative estimate takes $x_{d+1},y_{d+1}\ge\min w_i$, and $\|\tilde{\mathbf{p}}(t)\|\le\max_j\|\mathbf{r}_j\|$. Substituting yields the stated bound.
\end{proof}

\begin{Theorem}
Let a rational B\'ezier curve of degree \(n\) be given by control points \(\mathbf{r}_i \in \mathbb{R}^d\) and positive weights \(w_i > 0\) for \(i = 0, \dots, n\). Define
\[
C = \max_{0 < i < n} \left( \frac{\|\Delta_2(\omega\mathbf{r})_i\|}{\displaystyle\min_{0\le k\le n} \omega_k} + \frac{|\Delta_2 \omega_i| \,\displaystyle\max_{0\le j\le n} \|\mathbf{r}_j\|}{\displaystyle\min_{0\le k\le n} \omega_k} \right).
\]
For a given tolerance \(\epsilon > 0\), set \(T = \epsilon / C\). Then the smallest integer \(m > n\) such that the distance between the curve and its degree‑elevated control polygon does not exceed \(\epsilon\) can be bounded by
\[
m = \min\{ m_{\text{even}},\, m_{\text{odd}} \},
\]
where
\[
m_{\text{even}} = 2\left\lceil \frac{A+1}{2} \right\rceil,\qquad
A = \frac{n(n-1)}{8T},
\]
and
\[
m_{\text{odd}} = 2\left\lceil \frac{\max\{n+1,\, m^*\} - 1}{2} \right\rceil + 1,\qquad
m^* = \frac{n(n-1) + \sqrt{n^2(n-1)^2 + 32T\, n(n-1)}}{16T}.
\]
\end{Theorem}

\begin{proof}
From the general bound derived in the previous sections, the distance between the original curve and its control polygon after degree elevation to degree \(m\) is bounded by
\[
\frac{n(n-1)}{m(m-1)}\,\mathbb{N}_{\infty}(m) \cdot C,
\]
where \(\mathbb{N}_{\infty}(m) = \dfrac{\lfloor m/2\rfloor\lceil m/2\rceil}{2m}\).  
We require this quantity to be no larger than \(\epsilon\), i.e.
\[
\frac{n(n-1)}{m(m-1)}\,\mathbb{N}_{\infty}(m) \le T.
\]

The constant \(\mathbb{N}_{\infty}(m)\) has a simple form depending on the parity of \(m\):
\[
\mathbb{N}_{\infty}(m) = 
\begin{cases}
\dfrac{m}{8}, & \text{if } m \text{ is even},\\[6pt]
\dfrac{m^2-1}{8m}, & \text{if } m \text{ is odd}.
\end{cases}
\]

For even \(m\),
substituting \(\mathbb{N}_{\infty}(m)=m/8\) gives
\[
\frac{n(n-1)}{8(m-1)} \le T \quad\Longrightarrow\quad m-1 \ge \frac{n(n-1)}{8T} = A.
\]
Hence the smallest even integer satisfying \(m \ge A+1\) is
\[
m_{\text{even}} = 2\left\lceil\frac{A+1}{2}\right\rceil.
\]

For odd \(m\) we have
\[
\frac{n(n-1)(m+1)}{8m^2} \le T
\;\Longrightarrow\; 8T m^2 - n(n-1)m - n(n-1) \ge 0.
\]
The quadratic equation \(8T m^2 - n(n-1)m - n(n-1)=0\) has the positive root
\[
m^* = \frac{n(n-1) + \sqrt{n^2(n-1)^2 + 32T\, n(n-1)}}{16T}.
\]
Because the quadratic coefficient is positive, the inequality holds for all \(m > m^*\).  
Additionally we need \(m > n\); therefore the smallest odd integer that is at least \(\max\{n+1,\, m^*\}\) is
\[
m_{\text{odd}} = 2\left\lceil\frac{\max\{n+1,\, m^*\} - 1}{2}\right\rceil + 1.
\]

Finally, the overall minimal degree \(m\) satisfying the requirement is the smaller of the two candidates, i.e. \(m = \min\{m_{\text{even}},\, m_{\text{odd}}\}\). This \(m\) is guaranteed to be larger than \(n\) by construction.
\end{proof}

\paragraph*{Remark}
In the unlikely event that both candidates are \(\le n\), one may take \(m = n+1\) and directly verify the inequality; it usually holds because the bound is mild for small increases in degree. In practical applications, a threshold $E$ is usually preset. When the degree elevation number $m - n \ge E $,  set $ e = E $.

\label{sec:min_degree}

\section{Numerical Validation}
\label{Numerical Validation}

In this section, we provide two more examples to illustrate the effectiveness of our method. The results of the examples were obtained by Matlab(R) running on a PC with Intel(R) Core i7 CPU (2.2 GHz) with 16 GB of RAM. The detailed algorithm is described below.
\begin{algorithm}[H]
\caption{Algorithm for estimating the derivative bound of a rational Bézier curve}
\label{alg:derivative_bound}
\begin{algorithmic}[1]
\Require Control points $\mathbf{b}_i \in \mathbb{R}^2$ and weights $w_i>0$ for a rational Bézier curve of degree $n$, tolerance $\varepsilon$
\Ensure Upper bound $B$ for $max\|\mathbf{r}'\|$
\State Compute derivative curve coefficients $(\mathbf{Q}_i, \omega_i)_{i=0}^{2n}$ using Theorem~1
\State Determine the degree elevation steps $e$ according to Theorem~7
\State Compute the upper bound of $\|\mathbf{r}'(t)\|$ based on Theorem~2
\State Compare with the exact value and the conjectured value of the first derivative
\State \Return $B$
\end{algorithmic}
\end{algorithm}

\begin{example}
  For rational Bézier curves spanning degrees $n = 2$ to $10$, we employed the following test configuration
   \begin{itemize}
\item \textbf{Weights:} $ \qquad \ \  \omega_i=\alpha^{|i-m|}, \  m=\lfloor n/2\rfloor,\ \alpha = 1.1.$
\item \textbf{Control points:} $\mathbf r_i = \bigl(i, \ 0.05\sin(\pi i/n)\bigr),$ for $i = 0, \ldots, n.$
\item \textbf{ tolerance:} $\varepsilon=0.01.$
\end{itemize}

\begin{table}[htbp]
\centering
\caption{Comparison with corrected Bounds formula}
\label{tab:result1}
\begin{tabular}{c c c c | c c c}
\hline
$n$ & \makecell{Max first\\derivative} & \makecell{Corresponding\\ $t$} & \makecell{Conjectured\\upper bounds}  & \makecell{Our method} & \makecell{Running times\\(sec)}& \makecell{degree elevation \\ number ($e$)} \\ 
\hline
2 & 2.095238 & 0.500000 & 2.202748 & 2.103801 & 0.019703 & 32 \\ 
3 & 3.186405 & 0.364670 & 3.303092 & 3.190757 & 0.000339 & 147 \\ 
4 & 4.292683 & 0.500000 & 4.402749 & 4.298076 & 0.000429 & 230 \\ 
5 & 5.415446 & 0.417260 & 5.502375 & 5.419033 & 0.000810 & 558 \\ 
6 & 6.556553 & 0.500000 & 6.602062 & 6.570380 & 0.000421 & 215 \\ 
7 & 7.706513 & 0.439140 & 7.701812 & 7.718378 & 0.000621 & 357 \\ 
8 & 8.875039 & 0.500000 & 8.801611 & 8.883881 & 0.001293 & 665 \\ 
9 & 10.049171 & 0.451520 & 9.901447 & 10.058799 & 0.001488 & 800 \\ 
10 & 11.241522 & 0.500000 & 11.001313 & 11.253941 & 0.001647 & 800 \\ 
\hline
\end{tabular}
\end{table}

\begin{figure}[htbp]
    \centering
    \begin{minipage}[b]{0.5\textwidth}
        \centering
        \includegraphics[width=\textwidth]{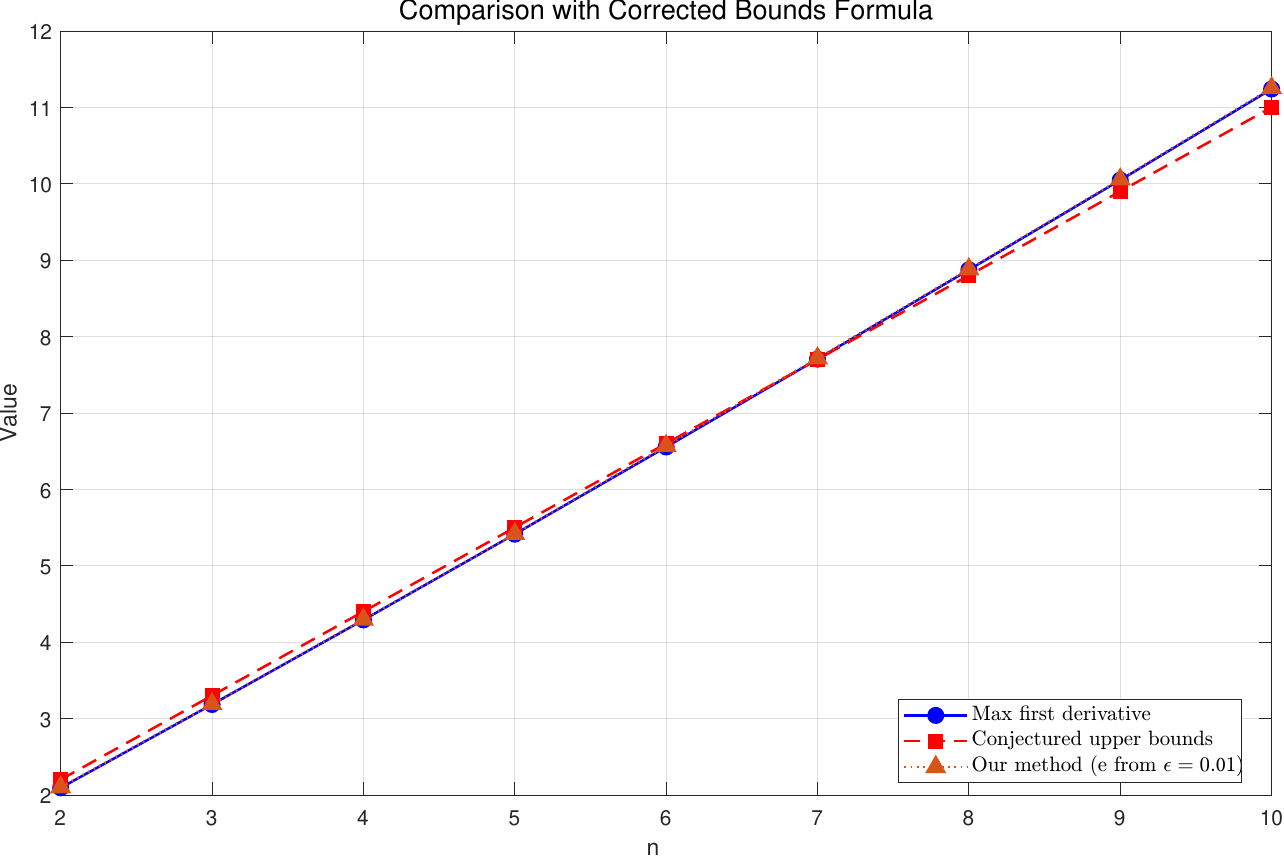}
        \vspace{2mm}
        \footnotesize (a) 
    \end{minipage}
    \hfill
    \begin{minipage}[b]{0.5\textwidth}
        \centering
        \includegraphics[width=\textwidth]{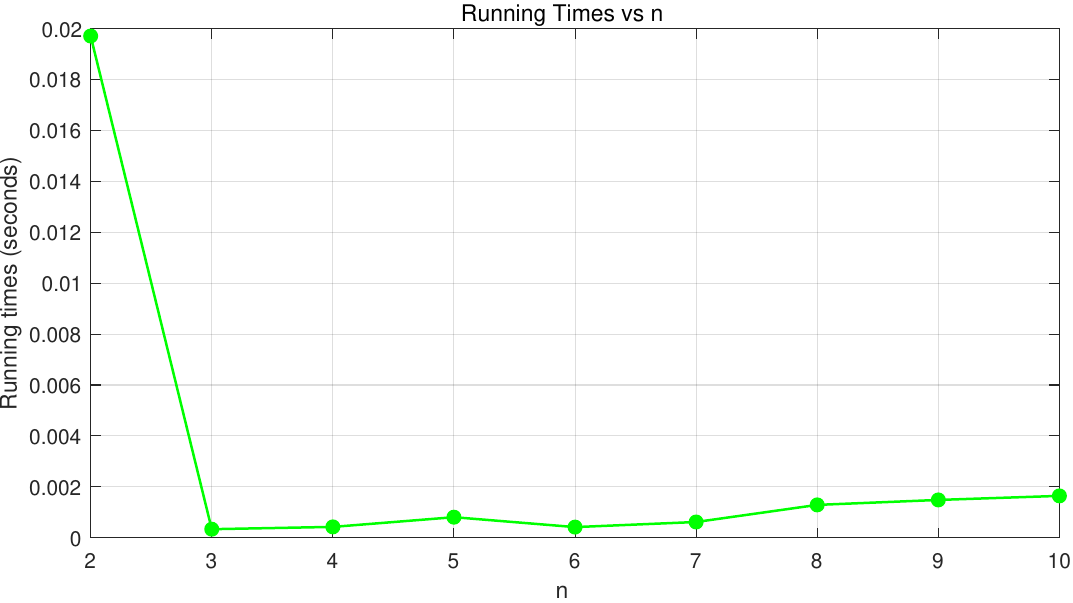}
        \vspace{2mm}
        \footnotesize (b) 
    \end{minipage}
    
    \begin{minipage}[b]{0.5\textwidth}
        \centering
        \includegraphics[width=\textwidth]{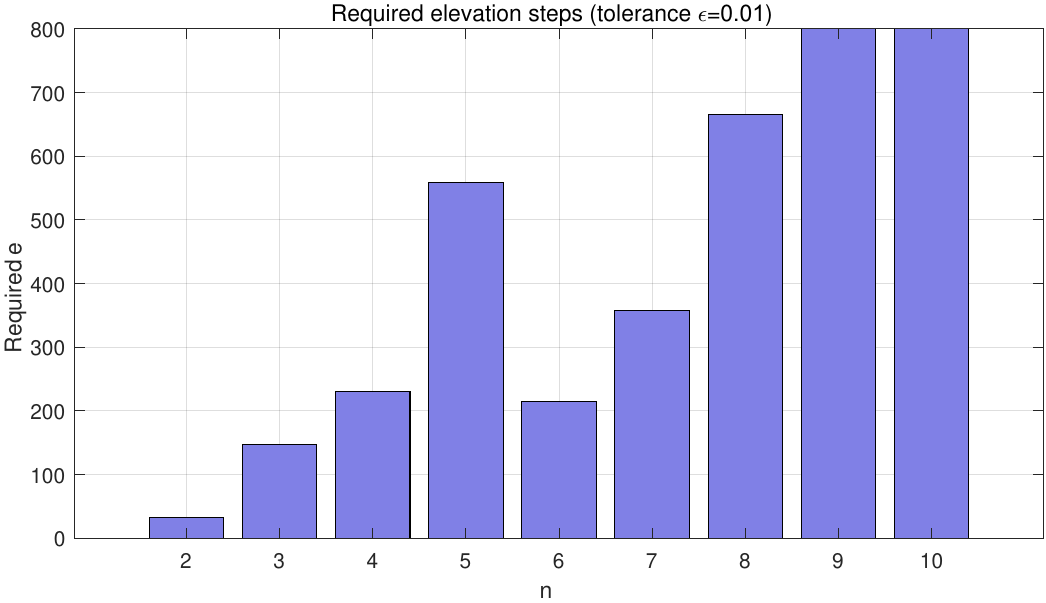}
        \vspace{2mm}
        \footnotesize (c) 
    \end{minipage}
    
    \vspace{3mm}
    \caption{Schematic of Example 2. (a) Comparison of derivative estimation methods. The plot shows the true maximum derivative, the conjectured upper bounds  and our computed bounds  for degrees $n = 2$ to $10$. For $n \geq 7 $, the true maximum exceeds the conjectured bound, while our method consistently provides valid upper bounds. (b) Computational performance of our method. In the tested range, the running time is nearly flat, indicating practical efficiency of the degree elevation approach. (c) The relationship between $\varepsilon=0.01$ and the degree elevation order $e$, where the threshold is set to $E=800$.}
    \label{fig2:comparison1}
\end{figure}

In this example, the degree elevation number is determined by $\epsilon$. Table~\ref{tab:result1} compiles our experimental findings, including maximum derivative values, corresponding parameter locations, conjectured bounds, our computed bounds,  computational times, and  degree elevation number. Figures~\ref{fig2:comparison1}  offer visual comparisons of different methodological approaches and their computational characteristics.
\end{example}

\begin{example}
  For rational Bézier curves spanning degrees $n = 2$ to $15$, we employed the following test configuration:

\begin{itemize}
\item \textbf{Weights:} $\omega_i = 2^{-i}$ for $i = 0, \ldots, n-1$, with $\omega_n = 2^{-(n-2)}$.
\item \textbf{Control points:} $\mathbf{r}_i = (i, 0) \in \mathbb{R}^2$ for $i = 0, \ldots, n$.
\item  \textbf{degree elevation number:} $e=800$.
\end{itemize}

\begin{table}[h!]
\caption{Comparison with corrected Bounds formula}
\label{tab:result2}
\centering
\begin{tabular}{c c c c | c c}
\hline
$n$ & \makecell{Max first\\derivative} & \makecell{Corresponding\\ $t$} & \makecell{Conjectured\\upper bounds} & \makecell{Our method\\(e=800)} & \makecell{Running times\\(sec)} \\ 
\hline
2 & 2.666667 & 0.500000 & 4.000000 & 2.669989 & 0.017191 \\ 
3 & 4.466444 & 0.656240 & 6.000000 & 4.475724 & 0.001169 \\ 
4 & 6.464102 & 0.732050 & 8.000000 & 6.482606 & 0.000907 \\ 
5 & 8.571204 & 0.778670 & 10.000000 & 8.602066 & 0.000988 \\ 
6 & 10.748531 & 0.810750 & 12.000000 & 10.794831 & 0.001028 \\ 
7 & 12.974278 & 0.834350 & 14.000000 & 13.038869 & 0.001212 \\ 
8 & 15.235043 & 0.852530 & 16.000000 & 15.320609 & 0.001548 \\ 
9 & 17.522048 & 0.867000 & 18.000000 & 17.631514 & 0.001416 \\ 
10 & 19.829270 & 0.878810 & 20.000000 & 19.965057 & 0.001603 \\ 
11 & 22.152423 & 0.888650 & 22.000000 & 22.317326 & 0.001833 \\ 
12 & 24.488370 & 0.896970 & 24.000000 & 24.684655 & 0.002560 \\ 
13 & 26.834755 & 0.904120 & 26.000000 & 27.064674 & 0.002072 \\ 
14 & 29.189773 & 0.910320 & 28.000000 & 29.456223 & 0.002336 \\ 
15 & 31.552017 & 0.915750 & 30.000000 & 31.857283 & 0.002630 \\ 
\hline
\end{tabular}

\end{table}

\begin{figure}[htbp]
    \centering
    \begin{minipage}[b]{0.5\textwidth}
        \centering
        \includegraphics[width=\textwidth]{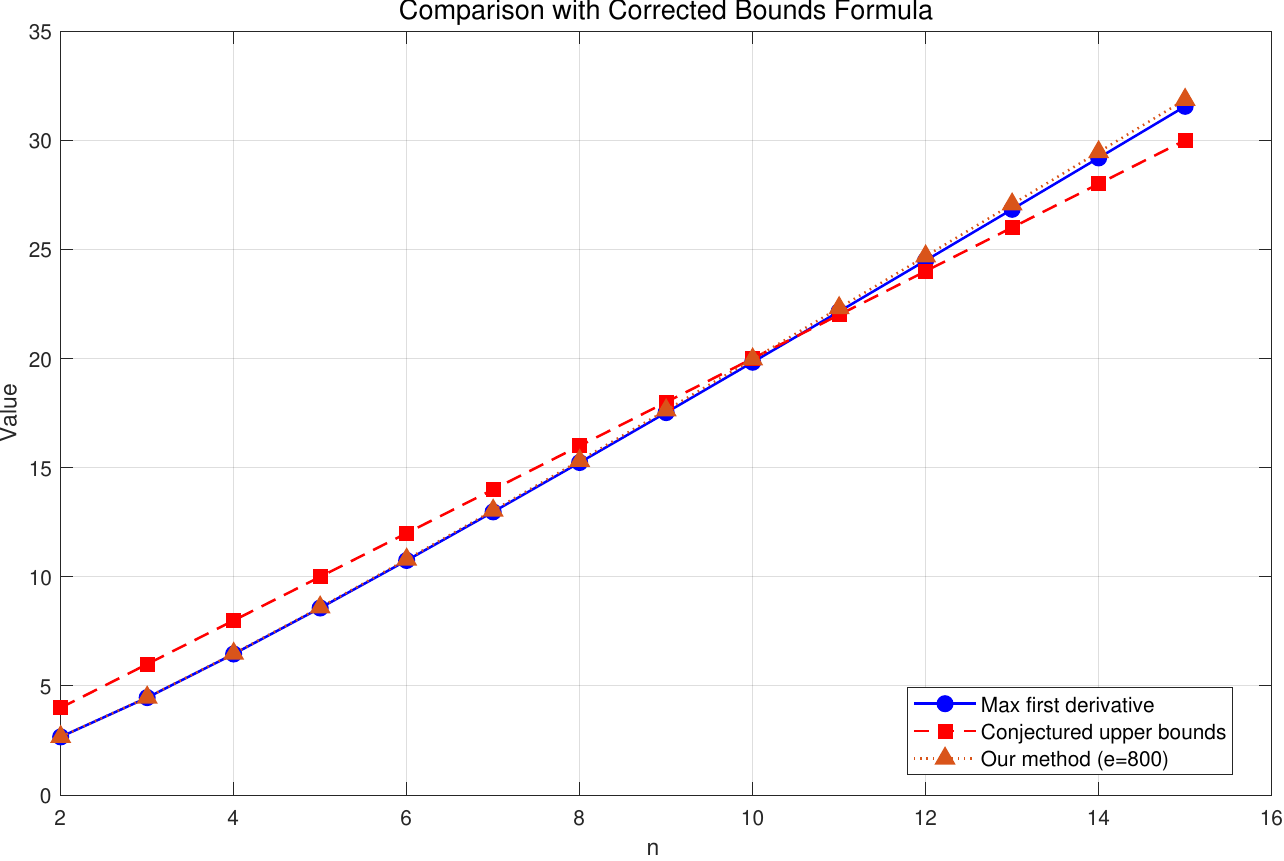}
        \vspace{2mm}
        \footnotesize (a) 
    \end{minipage}
    \hfill
    \begin{minipage}[b]{0.5\textwidth}
        \centering
        \includegraphics[width=\textwidth]{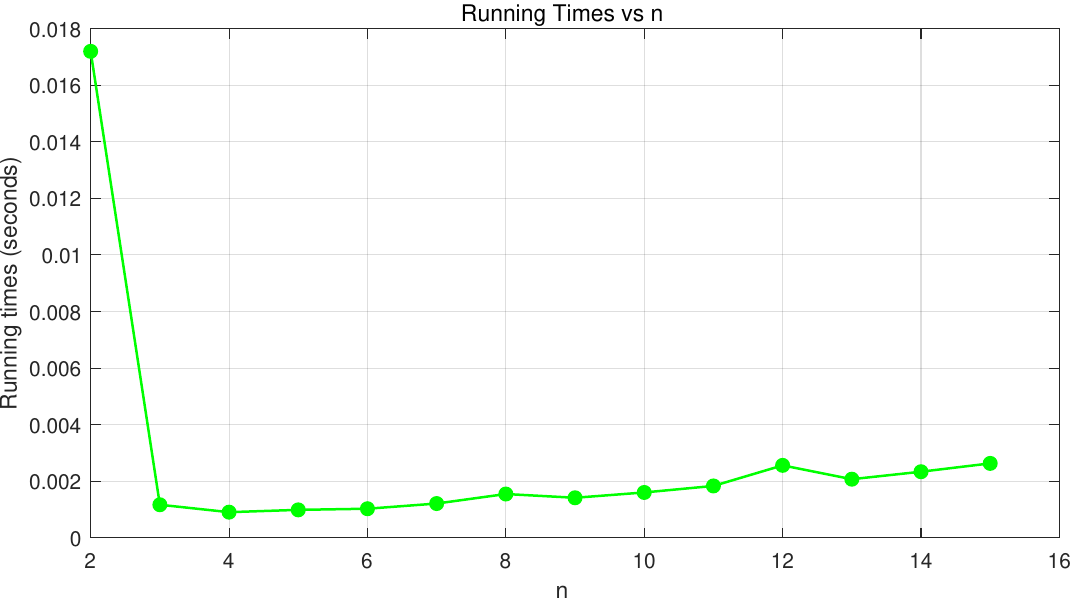}
        \vspace{2mm}
        \footnotesize (b) 
    \end{minipage}
    
    \vspace{3mm}
    \caption{Schematic of Example 3. (a) Comparison of derivative estimation methods. The plot shows the true maximum derivative (circles), the conjectured upper bounds are dashed line, and our computed bounds are triangles  for degrees $n = 2$ to $15$. For $n \geq 11 $, the true maximum exceeds the conjectured bound, while our method consistently provides valid upper bounds. (b) Computational performance of our method. In the tested range, the running time is nearly flat, indicating practical efficiency of the degree elevation approach with fixed $e = 800$. }
    \label{fig3:comparison1}
\end{figure}

In this example, for simplicity, we directly take the degree elevation number $e$ as a constant. This configuration yields $M = 2$ and $\max_{0 \leq i \leq n-1}\norm{\mathbf{r}_{i+1} - \mathbf{r}_i} = 1$, thus the conjectured bound equals $2n$.
Table~\ref{tab:result2} compiles our experimental findings, including maximum derivative values, corresponding parameter locations, conjectured bounds, our computed bounds, and computational requirements. Figures~\ref{fig3:comparison1}  offer visual comparisons of different methodological approaches and their computational characteristics.

\end{example}

\section{Conclusion and Open Problems}

This investigation successfully refutes the conjecture concerning derivative bounds of rational Bézier curves through carefully constructed counterexamples. Our analysis brings to light several significant unresolved questions:

\begin{itemize}
\item What represents the maximum degree $n$ for which the conjecture maintains validity? In addition, one of the reviewers of the paper has raised the following question:
What are the conditions under which the conjecture still holds when $n>6$?  e.g., monotone weights, bounded ratio sequences, etc..
\item What underlying mathematical principles account for the conjecture's failure? Might this relate to non-uniform convergence properties of rational Bézier curves under weight modifications \cite{Shi2005}?
\item Can we establish more restrictive, mathematically rigorous bounds that accommodate the observed supremum behavior of $\norm{\mathbf{r}'(t)}$?
\item This paper presents a relationship between the elevation number \(e\) and \(\varepsilon\), but it is too loose and not very effective in practical applications. Therefore, how to determine a smaller $e$ is also the goal of our future research.
\end{itemize}

\bibliography{mybibfile}
\bibliographystyle{unsrt}
 
\end{document}